\documentclass[a4paper,11pt]{article}
\usepackage{titlesec}
\usepackage{cfr-lm}
\usepackage{stmaryrd}
\usepackage{cite}
\usepackage[all]{xy}
\usepackage[dvips]{graphicx}
\usepackage{graphpap}
\usepackage{ifthen}
\usepackage{enumerate}
\usepackage{amssymb}
\usepackage{mathrsfs}
\usepackage[mathscr]{eucal}
\usepackage[errorshow]{tracefnt}
\usepackage{fancybox}
\usepackage{amsmath}
\usepackage{amssymb}
\usepackage{multirow}
\usepackage{array}
\usepackage{mathtools}
\usepackage{soul}
\usepackage{pict2e}
\usepackage{fonttable}
\DeclareMathAlphabet{\mathpzc}{OT1}{pzc}{m}{it}

\DeclareFontFamily{U}{matha}{\hyphenchar\font45}
\DeclareFontShape{U}{matha}{m}{n}{
      <5> <6> <7> <8> <9> <10> gen * matha
      <10.95> matha10 <12> <14.4> <17.28> <20.74> <24.88> matha12
      }{}
\DeclareSymbolFont{matha}{U}{matha}{m}{n}

\DeclareMathSymbol{\pm}            {2}{matha}{"08}
\DeclareMathSymbol{\mp}            {2}{matha}{"09}
\DeclareMathSymbol{\varleftarrow}{3}{matha}{"D0}
\DeclareMathSymbol{\varrightarrow}{3}{matha}{"D1}
\DeclareMathSymbol{\vee}           {2}{matha}{"5F}
\DeclareMathSymbol{\wedge}         {2}{matha}{"5E}
\DeclareMathSymbol{\leq}         {3}{matha}{"A4}
\DeclareMathSymbol{\geq}         {3}{matha}{"A5}
\DeclareMathSymbol{\in}            {3}{matha}{"50}
\DeclareMathSymbol{\owns}          {3}{matha}{"51}

\usepackage{bbold}

\usepackage{tikz}

\makeatletter
\DeclareRobustCommand{\Lcorner}{\mathbin{\mspace{1mu}\text{\L@corner}\mspace{1mu}}}
\newcommand{\L@corner}{%
  \setlength{\unitlength}{\fontcharht\font`T}%
  \begin{picture}(0.8,0)
  \roundcap
  \Line(0.1,0.95)(0.1,0.05)
  \Line(0.1,0.05)(0.7,0.05)
  \end{picture}%
}
\makeatother

\makeatletter
\DeclareRobustCommand{\Tri}{\mathbin{\mspace{1mu}\text{\L@corneer}\mspace{1mu}}}
\newcommand{\L@corneer}{%
  \setlength{\unitlength}{\fontcharht\font`T}%
  \begin{picture}(0.8,0)
  \roundcap
  \Line(0.1,0.052)(1.0,0.052)
  \end{picture}%
}
\makeatother

\newcommand{\Loc}{{\scriptscriptstyle\mkern 1mu\Lcorner}}
\newcommand{\LLoc}{{\scriptscriptstyle\mkern 1mu\mathbb{L}}}
\newcommand{\loc}{{\scriptscriptstyle\mkern 1mu\ell}}
\newcommand{\tri}{{\wedge\scriptscriptstyle\mkern -13.9mu\Tri}}

\makeatletter
\newcommand{\thickhline}{%
    \noalign {\ifnum 0=`}\fi \hrule height 1pt
    \futurelet \reserved@a \@xhline
}
\newcolumntype{'}{@{\hskip\tabcolsep\vrule width 1pt\hskip\tabcolsep}}
\makeatother
\newcolumntype{"}{@{\hskip\tabcolsep\vrule width 1.5pt\hskip\tabcolsep}}
\makeatother

\newcommand{\scr}{\mathscr}
\newcommand{\scrG}{\mathscr{G}}
\newcommand{\scrB}{\mathscr{B}}

\newcommand{\mangd}{M}

\newcommand{\crt}{{\!\vee\!}}

\def\boxit#1{\vbox{\hrule\hbox{\vrule\kern3pt
             \vbox{\kern3pt#1\kern3pt}\kern3pt\vrule}\hrule}}

\newcommand{\beq}{\begin{equation}}
\newcommand{\beqn}{\begin{equation*}}
\newcommand{\eeq}{\end{equation}}
\newcommand{\eeqn}{\end{equation*}}
\newcommand{\beqa}{\begin{eqnarray}}
\newcommand{\beqan}{\begin{eqnarray*}}
\newcommand{\eeqa}{\end{eqnarray}}
\newcommand{\eeqan}{\end{eqnarray*}}
\newcommand{\bdm}{\begin{displaymath}}
\newcommand{\edm}{\end{displaymath}}

\newcommand{\ba}{\begin{array}}
\newcommand{\ea}{\end{array}}

\newcommand\nn{\nonumber}

\newcommand\benu{\begin{enumerate}}
\newcommand\eenu{\end{enumerate}}
\newcommand\bit{\begin{itemize}}
\newcommand\eit{\end{itemize}}

\newtheorem{theorem}{Theorem}[section]
\newtheorem{lemma}[theorem]{Lemma}
\newtheorem{cor}[theorem]{Corollary}
\newtheorem{prop}[theorem]{Proposition}

\def\Pf{\noindent \textbf{Proof. }}

\def\der'{\mathfrak{der}'\,}
\def\der{\mathfrak{der}\,}
\def\str'{\mathfrak{str}'\,}
\def\str{\mathfrak{str}\,}

\def\qed{\hspace{\stretch{1}} $\Box$ \\
\noindent}

\newcommand{\de}{\delta}

\newcommand{\dlb}{\llbracket}%{\ensuremath{[\![}}
\newcommand{\drb}{\rrbracket}%{\ensuremath{]\!]}}

\newcommand{\blb}%{\ensuremath
{\text{$\llbracket$\hspace{-4pt}\scalebox{0.99}{$|$}\hspace{-2.58pt}\scalebox{0.99}{$|$}\hspace{-2.58pt}\scalebox{0.99}{$|$}}}
\newcommand{\brb}%{\ensuremath{
{\text{$\rrbracket$\hspace{-4pt}\scalebox{0.99}{$|$}\hspace{-2.58pt}\scalebox{0.99}{$|$}\hspace{-2.58pt}\scalebox{0.99}{$|$}}}

\def\fg{{\mathfrak g}}

\def\*{\partial}

\numberwithin{equation}{section}

\RequirePackage{hyperref} 

\begin{document}

\frenchspacing

%\allsectionsfont{\sffamily}

\titleformat{\section}
{\normalfont\bfseries\large\sffamily}{\textsf{{\thesection}}}{1em}{}
\titlelabel{\textsf{{\thetitle}}\quad}

\titleformat{\subsection}
{\normalfont\bfseries\sffamily}{\textsf{{\thesubsection}}}{1em}{}
\titlelabel{\textsf{{\thetitle}}\quad}

\titleformat{\subsubsection}
{\normalfont\bfseries\sffamily}{\textsf{{\thesubsubsection}}}{1em}{}
\titlelabel{\textsf{{\thetitle}}\quad}

\includegraphics[height=2cm]{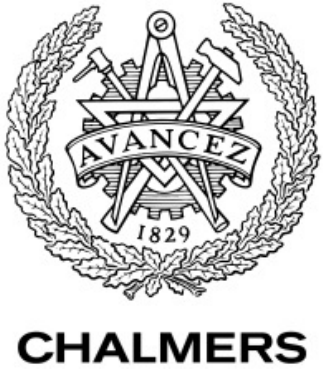}
\hspace{2mm}
\includegraphics[height=1.85cm]{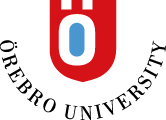}

%\vskip-10pt
%\hfill {\tt \today} \\
%\vskip-10pt
%\hfill {\tt \jobname} \\

\vspace*{1cm}

\vspace*{1.5cm}

\begin{center}
\noindent
{\LARGE {\sf \textbf{Tensor hierarchy algebras and restricted associativity}}}\\
\vspace{.3cm}

\renewcommand{\thefootnote}{\fnsymbol{footnote}}

\vskip 1truecm

\noindent
{\large {\sf \textbf{Martin Cederwall${}^{\,\mathsf 1}$ and Jakob Palmkvist${}^{\,\mathsf 2}$}
}}\\
\vskip .5truecm
        ${}^{\mathsf 1\,}${\it{Department of Physics\\ Chalmers University of Technology\\ SE-412 96 Göteborg, Sweden}\\[3mm]}
        {\tt martin.cederwall@chalmers.se}\\
\noindent
\vskip .5truecm
        ${}^{\mathsf 2\,}${\it{School of Science and Technology\\\"Orebro University\\ SE-701 82 \"Orebro, Sweden}\\[3mm]}
        {\tt jakob.palmkvist@oru.se} \\
\end{center}

\vskip 1cm

\centerline{\sf \textbf{
Abstract}}
\vskip .2cm

We study local algebras, which are structures similar to $\mathbb{Z}$-graded algebras
concentrated in degrees $-1,0,1$, but without a product
defined for pairs of elements at the same degree $\pm1$. 
To any triple consisting of
a Kac--Moody algebra $\fg$ with an
invertible and symmetrisable Cartan matrix, a dominant integral weight of $\fg$ and
an invariant symmetric bilinear form
on $\fg$, we associate a local algebra
satisfying a restricted version of associativity.
From it, we derive a local Lie superalgebra 
by a commutator construction. Under certain conditions, we identify generators which we show satisfy the relations of
the tensor hierarchy algebra $W$ previously defined from the same data.
The result suggests that an underlying structure satisfying such a restricted associativity
may be useful in applications of tensor hierarchy algebras to extended geometry.

\noindent

\newpage

\pagestyle{plain}

\tableofcontents

\section{Introduction}

The concept of local Lie algebras have played an important role in the
 classification of simple irreducible $\mathbb{Z}$-graded Lie algebras \cite{Kac68}
 (and thus to the development of Kac--Moody algebras)
by providing a ``seed'' at degrees $-1,0,1$ in the construction.
The concept can obviously be generalised from $\mathbb{Z}$-graded Lie algebras
to general $\mathbb{Z}$-graded algebras. However, it seems that such
``local algebras'' 
have not been studied much in 
cases other than those where the $\mathbb{Z}$-graded algebra is a Lie algebra or a Lie superalgebra \cite{Kac77B}.
Still in the context of Lie (super)algebras, it might for example be interesting to consider the commutator 
in a ``local associative algebra''.
It then turns out that the associative law is relevant only when at least one of the three involved elements
has degree zero. 
In the present paper, we introduce the concept of {\it focal associativity}
for local algebras where the associative law is restricted to these cases.
We will show that such a structure can be seen as
underlying
{\it tensor hierarchy algebras}, which are
infinite-dimensional generalisations of Cartan-type Lie superalgebras \cite{Palmkvist:2013vya,Carbone:2018xqq,Cederwall:2021ymp}. 
Tensor hierarchy algebras, originally used in the context of gauged supergravity \cite{Greitz:2013pua},
have
proven very
useful in the framework of extended geometry, where diffeomorphisms are unified with gauge transformations
in supergravity theories \cite{Bossard:2017wxl,Bossard:2019ksx,Cederwall:2019qnw,Cederwall:2019bai,Bossard:2021ebg,Cederwall:2021xqi}.

The paper is organised as follows. 
In section 2 we introduce the concept of local algebras, generalising the concept of
local Lie algebras introduced by Kac \cite{Kac68}, which we also specialise to
{\it contragredient} local Lie superalgebras.
We will show how any contragredient local Lie superalgebra $\scr G^\LLoc$
gives rise to a focally associative local superalgebra $\scr G^\loc$, which in turns gives back a different 
local Lie superalgebra $\scr G^\Loc$ with the commutator in $\scr G^\loc$ as the bracket.
In section 3 we show how a contragredient local Lie superalgebra $\scr B^\LLoc$ can be defined from a triple $(\fg,\lambda,\kappa)$,
where $\fg$ is a symmetrisable Kac--Moody algebra, $\lambda$ is a 
dominant integral weight of $\fg$ and $\kappa$
is an invariant symmetric bilinear form on $\fg$. This contragredient local Lie superalgebra $\scr B^\LLoc$
is the local part of a contragredient Lie superalgebra $\scr B$, which is also a Borcherds--Kac--Moody superalgebra \cite{Kac77B}.
In section~4 we 
recall the definition by generators
and relations of a
tensor hierarchy algebra $W$ from the same data $(\fg,\lambda,\kappa)$, under some further conditions 
\cite{Carbone:2018xqq}. We then apply the
construction in section 1 to the contragredient Lie superalgebra $\scr B^\LLoc$ defined in section 2. We identify the generators of $W$
with elements in $\scr B^\Loc$ and show that they generate a subalgebra where the defining relations of $W$ are satisfied up to an ideal 
intersecting the degree-zero subspace trivially. Throughout the paper, the base field $\mathbb{K}$ is algebraically closed and of characteristic zero.

\section{Local algebras}

We start by recalling that a {\it$\mathbb{Z}$-graded algebra} is a $\mathbb{Z}$-graded vector space $U=\bigoplus_{k \in \mathbb{Z}}U_k$
together with a degree-preserving map $U \otimes U \to U$, where the $\mathbb{Z}$-grading on $U \otimes U$
is given by $(U \otimes U)_k = \bigoplus_{i+j=k}U_i \otimes U_j$. Similarly, we define a {\it local algebra} as
a $\mathbb{Z}$-graded vector space $U=U_{-1}\oplus U_0 \oplus U_1$ together with a degree-preserving map $\bigoplus_{k=-1,0,1}(U \otimes U)_k \to U$.
The image of a pair $(u,v)$
is generally called {\it product} (as well as the map itself) and denoted $uv$, 
but in the Lie cases below, it will be called {\it bracket} and denoted $[u,v]$ or $\dlb u,v \drb$.
Note that a local algebra is actually not an algebra since the product is not defined for any pair of elements.

In a $\mathbb{Z}$-graded or local {\it superalgebra}, the product is also degree-preserving with respect to an additional
$\mathbb{Z}_2$-grading, $U=U_{(0)}\oplus U_{(1)}$.
The $\mathbb{Z}$-grading
is {\it consistent}
if $U_i \subseteq U_{(j)}$ whenever $i \equiv j$ (mod $2$).
In powers of $-1$, we will simplify the notation and write,
for example, $(-1)^{uv}$ for homogeneous elements $u,v$, where the exponent is actually the product of
their  $\mathbb{Z}_2$-degrees. We will also use subscripts to denote $\mathbb{Z}$-degrees of
homogeneous components of elements, for example, $u=\bigoplus_{k\in\mathbb{Z}}u_k$ in a $\mathbb{Z}$-graded superalgebra, and 
$u=\bigoplus_{k=-1,0,1}u_k$ in a local superalgebra.

Clearly, any $\mathbb{Z}$-graded algebra $U$ gives rise to a local algebra by restricting the vector space to the subspace 
$U=U_{-1}\oplus U_0 \oplus U_1$ and the domain of the product to $\bigoplus_{k=-1,0,1}(U \otimes U)_k$. 
This local algebra is called the {\it local part}
of the $\mathbb{Z}$-graded superalgebra $U$.

We say that a local algebra is {\it focally associative} if the degree-zero subspace associates with any element, that is, if
the identity
$(u_iv_j)w_k = u_i (v_jw_k)$ 
holds whenever all involved products are defined and
at least one of the three indices $i,j,k$ is zero. 
Thus the following 
13 identities are satisfied for any $u,v,w$ in a focally associative
local algebra,
%\vspace{.3cm}
\begin{align}
(u_{0}v_0)w_0&=u_{0}(v_0w_0)\,,\label{vanligass}\\[8pt]
%\end{align}
%\vspace{-.6cm}
%\begin{align}
(u_{\pm1}v_0)w_0&=u_{\pm1}(v_0w_0)\,,\label{vanster}\\
(u_{0}v_0)w_{\pm1}&=u_{0}(v_0w_{\pm1})\,,\label{hoger}\\
(u_{0}v_{\pm1})w_0&=u_{0}(v_{\pm1}w_0)\,,\label{mitten}\\[8pt]
%\end{align}
%\vspace{-.6cm}
%\begin{align}
(u_0v_{\pm1})w_{\mp1}&=u_0(v_{\pm1}w_{\mp1})\,,\label{vansterpm}\\
(u_{\pm1}v_{\mp1})w_0&=u_{\pm1}(v_{\mp1}w_0)\,,\label{hogerpm}\\
(u_{\pm1}v_{0})w_{\mp1}&=u_{\pm1}(v_{0}w_{\mp1})\,.\label{mittenpm}
\end{align}
If in addition the two identities
\begin{align} \label{extrass}
(u_{\pm1}v_{\mp1})w_{\pm1} = u_{\pm1}(v_{\mp1}w_{\pm1})
\end{align}
are satisfied for any $u,v,w$, then the local algebra is {\it associative}.

We define concepts like {\it subalgebras}
and {\it ideals} of local algebras in the same way as of algebras.
For any ideal $D$ of a local algebra
$U$, we also define the {\it quotient algebra} $U/D$ in the same way as for an ideal of an algebra. (Thus subalgebras and quotient algebras
of local algebras are local algebras as well.)
We say that the ideal $D$ is {\it peripheral} if $D=D_{-1} \oplus D_1$, where $D_{\pm1} \subseteq U_{\pm1}$.
The sum of all peripheral ideals is again a peripheral ideal, and therefore a unique maximal peripheral ideal.

Let $U$ be a local algebra and let $M$ be a subset of it. With
{\it the subalgebra of $U$ generated by $M$ modulo the maximal peripheral ideal}
we mean the quotient
algebra $V/D$, where $V$ is the subalgebra of $U$ generated by the subset
$M$, and $D$ is the maximal peripheral ideal of $V$.

A {\it local Lie superalgebra} (the logical ordering of words from our point of view here
would rather be {\it Lie local superalgebra}, but we stick to the established one) is a local superalgebra where the product is a 
{bracket} that satisfies the 
graded antisymmetry
\begin{align}
[x,y]=-(-1)^{xy}[y,x]
\end{align}
and the Jacobi identity
\begin{align}
[[x,y],z]&=[x,[y,z]]-(-1)^{xy}[y,[x,z]]
\end{align}
for any homogeneous elements such that
the involved brackets are defined. 
These two identities
can be broken down into the three plus five identities
\begin{align}
[x_0,y_0]&=-(-1)^{xy}[y_0,x_0]\,,\nn\\
[x_{\pm1},y_{\mp1}]&=-(-1)^{xy}[y_{\mp1},x_{\pm1}]\,,\\[8pt]
%\end{align}
%\vspace{-.8cm}
%\begin{align}
[[x_0,y_0],z_0]&=[x_0,[y_0,z_0]]-(-1)^{xy}[y_0,[x_0,z_0]]\,,\nn\\
[[x_0,y_0],z_{\pm1}]&=[x_0,[y_0,z_{\pm1}]]-(-1)^{xy}[y_0,[x_0,z_{\pm1}]]\,,\nn\\
[[x_{\pm1},y_{\mp1}],z_0]&=[x_{\pm1},[y_{\mp1},z_0]]-(-1)^{xy}[y_{\mp1},[x_{\pm1},z_0]]\,
\end{align}
for elements
that are homogeneous 
not only with respect to the $\mathbb{Z}_2$-grading, but
also with respect to the $\mathbb{Z}$-grading,
in the same way as the associative identity $(uv)w=u(vw)$ can be broken down into the 15 identities 
(\ref{vanligass})--(\ref{extrass}) above.

For any local superalgebra $\scr G^\loc$, we let $\scr G^\Loc$ be the superalgebra
which is the same vector space as $\scr G^\loc$, but 
with a different product which is
a bracket given by
the {\it commutator} $[x,y]=xy-(-1)^{xy}yx$ for homogeneous elements $x,y$.
The following proposition is an immediate consequence of the corresponding fundamental
statement for associative and Lie superalgebras, and straightforward to prove.

\begin{prop} \label{faimplieslie}
If $\scr G^\loc$ is focally associative, then $\scr G^\Loc$ is a local Lie superalgebra.
\end{prop}
\noindent
The reason why focal associativity is sufficient is that there is no Jacobi identity involving two elements
at degree $\pm1$ and one element at degree $\mp1$, since such an identity would involve the bracket of
the two elements at degree $\pm1$, which is not defined in a local Lie superalgebra.

We will now go in the opposite direction and associate a focally associative local algebra to 
a local Lie superalgebra satisfying some further conditions.

Let $\mathscr G^\LLoc=\mathscr G^\LLoc{}_{-1} \oplus \mathscr G^\LLoc{}_{0} \oplus \mathscr G^\LLoc{}_{1}$ be
a local Lie superalgebra with a bracket $\dlb-,-\drb$.
We say that $\mathscr G^\LLoc$ is {\it contragredient} if
there is an element $L \in \scr G_0$ such that
$\dlb L,x_k\drb=kx_k$ for all $x \in \scr G$ and a
bilinear map
\begin{align}
\mathscr G_{-1} \times \mathscr G_{1} \to \mathbb{K}\,, \qquad (x,y) \mapsto \langle x | y \rangle\,,
\end{align}
which 
is {\it invariant} and {\it homogeneous}. The conditions of invariance and homogeneity mean, respectively, that
$\langle \dlb x_{-1},y_0 \drb|z_1 \rangle = \langle x_{-1} | \dlb y_0,z_1\drb \rangle$
for all $x,y,z \in \scrG^\LLoc$,
and that $\langle x | y \rangle=0$ whenever $x$ and $y$ are homogeneous with different $\mathbb{Z}_2$-degrees.
It is convenient to also define a corresponding bilinear map
\begin{align}
\mathscr G_{1} \times \mathscr G_{-1} \to \mathbb{K}\,, \qquad (x,y) \mapsto \langle x | y \rangle = (-1)^{xy} \langle y | x \rangle\,.
\end{align}
by graded symmetry.

To any contragredient local Lie superalgebra $\scrG^\LLoc=\scrG^\LLoc{}_{-1}\oplus\scrG^\LLoc{}_0\oplus\scrG^\LLoc{}_1$, we associate a focally associative
local superalgebra $\scrG^\loc=\scrG^\loc{}_{-1}\oplus\scrG^\loc{}_0\oplus\scrG^\loc{}_1$ in the following way.
Let $\scrG^\loc{}_0$ be the universal enveloping algebra of $\scrG^\LLoc{}_0$, set $\scrG^\loc{}_{\pm1}=\scrG^\LLoc{}_{\pm1}\otimes \scrG^\loc{}_0$
and write $x \otimes 1 =x$ for any $x \in \scrG^\LLoc{}_{\pm1}$ 
(so that we consider $\scrG^\LLoc{}_{\pm1}$ as a subspace of $\scrG^\loc{}_{\pm1}$).
Accordingly, we consider $\scrG^\LLoc{}_{-1}\oplus\mathbb{K}\oplus\scrG^\LLoc{}_{1}$
as a subspace of $\scrG^\loc$. For $x$ and $y$ in this subspace, set
\begin{align}
x_{-1}y_{1} &= - a\dlb x_{-1},y_{1}\drb+b\langle x_{-1}|y_{1}\rangle L\,,\label{xyminusplus}\\
x_{1}y_{-1} &= a\dlb x_{1},y_{-1}\drb+b\langle x_{1}|y_{-1}\rangle L +c\langle x_{1}|y_{-1}\rangle\,\label{xyplusminus}
\end{align}
for some constants $a,b,c \in \mathbb{K}$,
and let $x_0y_{\pm1}=y_{\pm1}x_0$ be given by the action of $x_0\in\mathbb{K}$ that
$\scrG^\LLoc{}_{-1}\oplus\scrG^\LLoc{}_{1}$ is equipped with as a vector space over $\mathbb{K}$.

Note that $\scrG^\LLoc{}_{-1}\oplus\mathbb{K}\oplus\scrG^\LLoc{}_{1}$ is in general not a local algebra with respect to the product defined so far,
since the right hand sides of (\ref{xyminusplus}) and (\ref{xyplusminus})
in general do not belong to this subspace of $\scrG^\loc$. In order to achieve a local algebra,
we will now extend the product to the whole of $\scrG^\loc$.
First, we 
recursively define subspaces $(\scrG^\loc{}_0)^k$ of $\scrG^\loc{}_0$ for any integer $k\geq0$
by setting  $(\scrG^\loc{}_0)^{0}=\mathbb{K}$ and
letting $(\scrG^\loc{}_0)^{k+1}$ consist of all
elements $ux$ where $u \in (\scrG^\loc{}_0)^{k}$ and $x \in \scrG^\LLoc{}_0$. 
As the universal enveloping algebra of $\scrG^\LLoc{}_0$,
the algebra $\scrG^\loc{}_0$ is the sum of all such subspaces.
Second, we define the product on $\scrG^\loc$ recursively by
\begin{align}
x \big(y\otimes v\big)=(xy)v \label{xyv}
\end{align}
and
\begin{align} \label{master}
\big(x \otimes (uz)\big)\big(y\otimes v\big)&=(x\otimes u)\big(\dlb z,y\drb\otimes v\big)+(-1)^{yz}(x\otimes u)\big(y\otimes (zv)\big)\,,
\end{align}
where
\begin{align}
x &\in \mathbb{K} \oplus \scrG^\LLoc{}_{\pm1}\,, & y &\in \mathbb{K} \oplus \scrG^\LLoc{}_{\mp1}\,, & z&\in\scrG^\LLoc{}_0\,,
& u &\in (\scrG^\loc{}_0)^i\,,& v &\in (\scrG^\loc{}_0)^j
\label{legend}
\end{align}
for $i\geq1$ and where we set $\dlb z,y\drb=0$ if $y \in \mathbb{K}$. It is straightforward to check that the product is well defined.

By setting 
$x=1$ in (\ref{xyv}), we see that the tensor product symbol $\otimes$ is superfluous (and it will henceforth be omitted).
Also, it follows from the two equations that we obtain from (\ref{master}) by setting $x=y=u=1$ and $x=y=v=1$
that the product on $\scrG^\ell{}_0$ defined by this equation is the same as the one that this vector space is equipped with 
as the universal enveloping algebra of $\scrG^\LLoc{}_0$.
The product is thus associative on $\scrG^\ell{}_0$.

Let us compute the commutator $[x,y]=xy-(-1)^{yx}yx$ given by the product above
for elements in $\scr G^\LLoc \subseteq \scrG^\loc$. It is
equal to the original bracket in the following cases,
\begin{align} \label{kommutatorlika}
[x_0,y_0]&=\dlb x_0,y_0 \drb \,,& [x_0,y_{\pm1}]&=\dlb x_0,y_{\pm1}\drb \,, & [x_{\pm1},y_0]&=\dlb x_{\pm1},y_0\drb\,,
\end{align}
but not when $x \in \scrG^\LLoc{}_{\pm1}$ and $y \in \scrG^\LLoc{}_{\mp1}$. In this case we instead get
\begin{align}
[x_{-1},y_1]&=x_{-1}y_1 - (-1)^{xy} y_1x_{-1}\nn\\
&=- a\dlb x_{-1},y_{1}\drb+b\langle x_{-1}|y_{1}\rangle L\nn\\
&\quad\,-a(-1)^{xy}\dlb y_{1},x_{-1}\drb-b(-1)^{xy}\langle y_{1}|x_{-1}\rangle L -c(-1)^{xy}\langle y_{1}|x_{-1}\rangle\nn\\
&=-c(-1)^{xy}\langle y_1 | x_{-1}\rangle =-c\langle x_{-1}|y_{1}\rangle\,. \label{kommutatorolika}
\end{align}

We will now show that the local algebra $\scrG^\loc=\scrG^\loc{}_{-1}\oplus\scrG^\loc{}_0\oplus\scrG^\loc{}_1$ is indeed focally associative.
We already know that the identity (\ref{vanligass}) holds for all $u,v,w \in\scrG^\loc$ since $\scrG^\loc$ is the universal enveloping algebra of
$\scr G^\LLoc$ and thus associative. The identities (\ref{vanster}), (\ref{hoger}) and (\ref{mittenpm}) are consequences of the following proposition.

\begin{prop} \label{flyttaoverw}
The identity
\begin{align} \label{identiteten}
\big((xu)w\big)(yv)&=(xu)\big(w(yv)\big)\,,
\end{align}
where
\begin{align}
x &\in \mathbb{K} \oplus \scrG^\LLoc{}_{\pm1}\,, & y &\in \mathbb{K} \oplus \scrG^\LLoc{}_{\mp1}\,, & 
u &\in (\scrG^\loc{}_0)^i\,, & v &\in (\scrG^\loc{}_0)^j\,, & w &\in (\scrG^\loc{}_0)^k\,,
\label{legend2}
\end{align}
holds for all integers $i,j,k \geq 0$.
\end{prop}
\Pf
We will prove this 
by induction over $i+k\geq0$. The base cases are trivial. Suppose the identities hold for $i+k\leq p$ for some $p\geq 0$.
For $i=p$ we then have
\begin{align}
((xu)z)(yv)=(x(uz))(yv)=(xu)([z,y]v)+(-1)^{yz}(xu)(y(zv))=(xu)(z(yv))\,,
\end{align}
where $z \in \scr G^\LLoc{}_0$,
by the induction hypothesis in the first step, and (\ref{master}) in the other two. Thus the identity (\ref{identiteten}) holds for $k=1$ and $i=p$.
It is now straightforward to proceed by induction over $k$, and show that it holds for any $k\geq 1$ and $i+k=p+1$. It 
suffices to say that the 
idea in the induction step of this second induction is, as in (\ref{flyttaover}) below, to move one element at the time from one pair of parentheses to the other. The proposition then follows by the principle of induction.
\qed

\noindent
We now turn to the remaining parts (\ref{mitten})--(\ref{hogerpm}) of the focal associativity.

\begin{lemma} \label{lemmat}
The identities
\begin{align}
(uy)v&=u(yv)\,,\nn\\
(ux)(yv)&=u(xy)v\label{naturbarn}
\end{align}
hold for all variables as in (\ref{legend2}).
\end{lemma}
\Pf
We suppose that $x \in \scrG^\LLoc{}_{\pm1}$ and $y \in \scrG^\LLoc{}_{\mp1}$, since this is sufficient, and prove the lemma 
by induction over $i$. The base case $i=0$ is either trivial or given by (\ref{xyv}).
Suppose the identities hold for $i\leq p$ for some $p\geq 0$.
For $i=p$ we then have
\begin{align}
(uzy)v
&=(u[z,y])v+(-1)^{yz}(uyz)v\nn\\
&=(u[z,y])v+(-1)^{yz}(uy)(zv)\nn\\
&=u([z,y]v)+(-1)^{yz}u(y(zv))=(uz)(yv)\,,\label{indstepnaturbarn0}
\end{align}
where we have used (\ref{identiteten}) in the second step, the induction hypothesis in the third and
(\ref{master}) in the fourth. For $i=p$ we furthermore have
\begin{align}
(uzx)(yv)&=(u[z,x])(yv)+(-1)^{zx}(uxz)(yv)\nn\\
&=u([z,x]y)v+(-1)^{zx}(ux)(zyv)\nn\\
&=u([z,x]y)v+(-1)^{zx}(ux)([z,y]v)+(-1)^{zx+zy}(ux)(yzv)\nn\\
&=u([z,x]y)v+(-1)^{zx}u(x[z,y])v+(-1)^{zx+zy}u(xy)zv\nn\\
&=u\big([z,x]y+(-1)^{zx}x[z,y]+(-1)^{zx+zy}(xy)z\big)v \label{indstepnaturbarn}
\end{align}
using the induction hypothesis in the second and fourth step. If $(x,y)=(x_{-1},y_{1})$, the expression between $u$ and $v$
equals
\begin{align}
[z,x]y&+(-1)^{zx}x[z,y]+(-1)^{zx+zy}(xy)z   \nn\\
&=-a\dlb \dlb z,x\drb,y\drb+b\langle [z,x]|y \rangle L\nn\\
&\quad\,-a(-1)^{zx}\dlb x,\dlb z,y\drb\drb+b(-1)^{zx}\langle x | \dlb z,y \drb\rangle L\nn\\
&\quad\,-a(-1)^{z(x+y)}\dlb x,y\drb z+b(-1)^{z(x+z)}\langle x | y \rangle zL\nn\\
&=-a\dlb z,\dlb  x,y\drb\drb
-a(-1)^{z(x+y)}\dlb x,y\drb z
+b(-1)^{z(x+z)}\langle x | y \rangle zL\nn\\
&=-az\dlb x,y\drb +b\langle x | y \rangle zL=z(xy)
\end{align}
Thus $(uzx)(yv)=uz(xy)v$.
The case $(x,y)=(x_{1},y_{-1})$ is similar. We have thus shown that the identities (\ref{naturbarn})
hold when $i=p+1$ as well, and the lemma follows by the principle of induction.
\qed

\noindent
Note that all products of three elements written without parentheses in (\ref{indstepnaturbarn0}) and (\ref{indstepnaturbarn}) are well defined,
either by the induction hypothesis or by Proposition \ref{flyttaoverw}. 

\begin{prop} \label{remainingthree}
The identities
\begin{align}
(u(yw))v&=u((yw)v)\,,\nn\\
((ux)(yv))w&=(ux)((yv)w)\,,\nn\\
(w(ux))(yv)&=w((xu)(yv))\,, \label{remaining}
\end{align}
hold for all variables as in (\ref{legend2}).
\end{prop}
\Pf We have
\begin{align} 
(u(yw))v&=((uy)w)v=(uy)(wv)=u(y(wv))=u((yw)v)\,,\nn\\
((ux)(yv))w&=u(xy)vw=(ux)(y(vw))=(ux)((yv)w)\,,\nn\\
(w(ux))(yv)&=((wu)x)(yv)=wu(xy)v=w\big((ux)(yv)\big)\,, \label{flyttaover}
\end{align}
by Proposition \ref{flyttaoverw} and Lemma \ref{lemmat}.
\qed 

\begin{cor}
The local algebra $\scrG^\loc$ is focally associative.
\end{cor}
\Pf This follows directly from Propositions \ref{flyttaoverw} and \ref{remainingthree}, and the fact that any element in $\scrG^\ell{}_{\pm1}$
can be written as a sum of elements $ux$, where $u\in\scrG^\ell{}_0$ and $x\in\scrG^\LLoc{}_{\pm1}$, which is
easily shown by induction.
\qed

\noindent
We have shown that $\scrG^\loc$ is a focally associative local superalgebra, and thus, by Proposition \ref{faimplieslie},
it gives rise to a new local Lie superalgebra
where the bracket is given by the commutator in $\scrG^\loc$. We denote this local Lie superalgebra by $\scrG^\Loc$,
and the bracket in it by $[-,-]$, to be distinguished from the original bracket $\dlb-,-\drb$ on $\scrG^\LLoc$. This is particularly important when
one of the elements belong to $\scrG^\LLoc{}_1$ and the other to $\scrG^\LLoc{}_{-1}$ since both brackets are defined in this case,
but disagree according to (\ref{kommutatorolika}).

Note that it was only in the second part of the proof of Lemma \ref{lemmat} that we used the form of the product $xy$ as an element in
$\mathbb{K}\oplus\scrG^\LLoc{}_0$, and that the values of the constants $a,b,c$ did not matter.
In fact, we can always assume $a=1$ without loss of generality by redefining the bracket $\dlb-,-\drb$. Similarly, whenever
$b\neq0$ or $c\neq0$, we can assume $b=1$ or $c=1$ without loss of generality by redefining the invariant form $\langle-|-\rangle$.
We will assume that both constants $b$ and $c$ are nonzero, and furthermore that they are equal to each other, since this
condition turns out be important for the relation to the tensor hierarchy algebras
(more precisely, it is crucial in the proof of Lemma \ref{sistalemmat} below). Accordingly, we henceforth set $a=b=c=1$,
and we have
\begin{align}
x_{-1}y_{1} &= - \dlb x_{-1},y_{1}\drb+\langle x_{-1}|y_{1}\rangle L\,,
\nn\\
x_{1}y_{-1} &= \dlb x_{1},y_{-1}\drb+\langle x_{1}|y_{-1}\rangle L +\langle x_{1}|y_{-1}\rangle\,,\label{xyplusminusutanabc}
\end{align}
so that
$[x_{\pm1},y_{\mp1}]=\pm\langle x_{\pm1}|y_{\mp1}\rangle$.
For any contragredient local Lie superalgebra $\scr G^\LLoc$ we thus
let $\scr G^\loc$ be the focally associative local algebra constructed in the way above with
this choice of constants $a,b,c$, and $\scrG^\Loc$ the local Lie superalgebra obtained from $\scrG^\loc$ with the commutator
$[-,-]$ as the bracket, to be distinguished from the original one $\dlb-,-\drb$.
Note that $\scrG^\loc$ (and thus also $\scrG^\Loc$) is in general infinite-dimensional even when
$\scr G^\LLoc$ is finite-dimensional.

\section{Contragredient Lie superalgebras}

Let $\fg$ be a Kac--Moody algebra of rank $r$ with an invertible and symmetrisable Cartan matrix $A$,
let $\lambda$ be a dominant integral weight of $\fg$ and
let $\kappa$ be an invariant symmetric bilinear form on $\fg$.
In this section we will associate 
a contragredient local Lie superalgebra $\scr B{}^\LLoc$ to the triple $(\fg,\lambda,\kappa)$,
from which we in turn can construct a focally associative local superalgebra
$\scr B^\loc$ and a local Lie superalgebra $\scr B^\Loc$ as above.

We recall that $\fg$ is generated by $3r$ elements $e_k,f_k,h_k$, where $k=1,2,\ldots,r$,
modulo the Chevalley--Serre relations \cite{Kac}.
We also recall that
the invariant symmetric bilinear form $\kappa$ on $\fg$ is unique up to an overall normalisation,
that it satisfies $\kappa(e_k,f_k)\neq 0$ for 
any $k=1,2,\ldots,r$, and that it induces a symmetric bilinear form on the vector space 
$\mathfrak{h}^\ast$ dual to the Cartan subalgebra $\mathfrak{h}$
(spanned by the generators $h_k$)
by the relation $(\alpha_i{}^\crt,\alpha_j{}^{\!\vee})=\kappa(h_i,h_j)$,
where the simple coroots are defined by $\alpha_k{}^\crt=\kappa(e_k,f_k)\alpha_k$.
It then follows that 
$(\alpha_k,\alpha_k)=2/{\kappa(e_k,f_k)}$
so that $\alpha_k{}^\crt=2\alpha_k/(\alpha_k,\alpha_k)$.
These well known results will be re-derived below for the contragredient Lie superalgebra $\scr B$
with Cartan matrix $B$ obtained by adding a row and column to the Cartan matrix $A$.

Let $\lambda_k =(\lambda,\alpha_k{}^\vee)$ be the 
Dynkin labels of the 
dominant integral weight $\lambda$, so that $\lambda_k\in\mathbb{Z}$ and $\lambda_k\geq0$ for any $k=1,2,\ldots,r$ (not all zero).
The Dynkin labels are the components of $\lambda$ in the basis of fundamental weights $\Lambda_k$, defined by $(\Lambda_i,\alpha_j{}^\crt)=\delta_{ij}$.
Let $\lambda^{\tri}$ be the weight with Dynkin labels $\lambda^{\tri}{}_k=\lambda_k/{\kappa(e_k,f_k)}$. 
We will be interested in cases where $\fg$ is finite and where $\lambda$ and $\kappa$ are such that $\lambda^\tri$ is 
a fundamental weight $\Lambda_k$ for which the corresponding Coxeter label (the component of the highest root
$\theta$ in the basis of simple roots) is equal to $1$.
We say that such a weight $\lambda^{\tri}$ is a {\it pseudo-minuscule} weight. The reason for choosing this term
(although it has been used in a different meaning \cite{Letzter200488})
is that the pseudo-minuscule weights coincide with the {\it minuscule} weights 
(highest weights of representations on which the Weyl group acts transitively \cite{Bourbaki}) for all $\fg$ other than $\fg=B_r$ and $\fg=C_r$.
Moreover, the isomorphism between the weight spaces of $B_r$ and $C_r$ given by 
transposing the Cartan matrix (or flipping the arrow in the Dynkin diagram) maps a minuscule weight of one algebra to a pseudo-minuscule weight of the other, and vice versa.
(This in fact holds for any Cartan matrix of a finite Kac--Moody algebra $\fg$, but for other $\fg$ it just says that the
minuscule and pseudo-minuscule weights coincide.) Below follows the 
complete list of pseudo-minuscule weights in the numbering of Bourbaki \cite{Bourbaki2}
(with some additional information about the corresponding highest weight representations).
There are no pseudo-minuscule weights of $E_8$, $F_4$ or $G_2$.
\begin{itemize}
\item $A_r:$ $\Lambda_1,\ldots,\Lambda_r$ 
\item $B_r:$ $\Lambda_1$ (vector representation)
\item $C_r:$ $\Lambda_r$
\item $D_r:$ $\Lambda_1,\Lambda_{r-1},\Lambda_r$ (vector and spinor representations)
\item $E_6:$ $\Lambda_1,\Lambda_6$ (27-dimensional)
\item $E_7:$ $\Lambda_7$ (56-dimensional)
\end{itemize}
In extended geometry with extended structure algebra $\fg$ and extended coordinate representation with highest weight
$\lambda$, it is precisely when $\lambda^\tri$
is a pseudo-minuscule weight that additional ``ancillary'' transformations are not needed for closure and covariance
of the generalised diffeomorphisms \cite{Cederwall:2017fjm}. (In \cite{Cederwall:2017fjm}, the normalisation was chosen such that
$\lambda=\lambda^\tri$, if possible. Accordingly, the conclusion there was that ancillary transformations are absent
precisely when $\lambda$ is a pseudo-minuscule weight. However, with a different normalisation they would presumably
be absent also when $\lambda$ is an integer multiple of a pseudo-minuscule weight.)

\begin{prop}
A necessary condition for {\rm$\lambda^\tri$} to be a pseudo-minuscule weight is that $(\lambda,\theta)=1$. If {\rm$\lambda^\tri$}
is a dominant integral weight, then this condition is also sufficient.
\end{prop}
\Pf If {\rm$\lambda^\tri$} is a pseudo-minuscule weight and $c_k$ are the component of $\theta$ in the basis of simple roots $\alpha_k$,
then
\begin{align}
1&=\sum_{k=1}^r \lambda^\tri{}_k c_k
=\sum_{k=1}^r \frac{(\alpha_k,\alpha_k)}2\lambda_k c_k=\sum_{i,j=1}^r \frac{(\alpha_j,\alpha_j)}2\lambda_i c_j\de_{ij}\nn\\&=
\sum_{i,j=1}^r \frac{(\alpha_j,\alpha_j)}2\lambda_i c_j(\Lambda_i,\alpha_j{}^\crt)=
\sum_{i,j=1}^r \lambda_i c_j(\Lambda_i,\alpha_j)=(\lambda,\theta)\,. \label{lambdatheta}
\end{align}
Conversely, if $(\lambda,\theta)=1$, then the same calculation shows that $\sum_{k=1}^r \lambda^\tri{}_k c_k=1$. If in addition
the Dynkin labels $\lambda^\tri{}_k$ are non-negative integers, then the only possibility is that all are zero except for
one of them which is equal to $1$, and that the corresponding Coxeter label $c_k$ is equal to $1$ too.
\qed

\noindent
Given the triple $(\fg,\lambda,\kappa)$, let $B$ be the square matrix of order $r+1$ with entries
\begin{align} \label{matrisdef}
B_{00}&=0\,,& B_{i0}&=-\lambda_i \,,& B_{0j}&=-\lambda^\tri{}_j=-\frac{\lambda_j}{\kappa(e_j,f_j)}\,& B_{ij}&=A_{ij}\,, %, &(e_0,f_0)&=1\,.
\end{align}
where $i,j=1,2,\ldots,r$.
Then $B$
is symmetrisable. We also assume that $\lambda$ is such that $B$ is invertible.

The {\it contragredient Lie superalgebra} $\scr B$ 
associated to the Cartan matrix $B$ is defined 
from a set of $3r$ generators $\mangd_{\scr B}=\{e_K,f_K,h_K | K=0,1,\ldots,r\}$, where
$e_0$ and $f_0$ are odd, whereas $h_0$ and $e_k,f_k,h_k$ are even for $k=1,2,\ldots,r$. Let $\tilde{\scr B}$ be the $\mathbb{Z}$-graded Lie superalgebra
generated by this set $M_{\scr B}$ modulo the relations 
\begin{align}
[h_I,e_J]&=B_{IJ}e_J\,, & [h_I,f_J]&=-B_{IJ}f_J\,, & [e_I,f_J]&=\delta_{IJ}h_J  
\end{align}
with the (non-consistent) $\mathbb{Z}$-grading where $e_K$ and $f_K$ have degree $1$ and $-1$, respectively, for any $K=0,1,\ldots,r$. 
Then $\scr B=\tilde{\scr B}/D$, where $D$ is the maximal graded ideal of $\tilde{\scr B}$ intersecting the local part of $\tilde{\scr B}$ trivially
\cite{Kac77B}.
Since $B$ here satisfies the conditions of a Cartan matrix of a {\it Borcherds--Kac--Moody algebra}, 
a generalisation \cite{Ray98} of the Gabber--Kac theorem \cite{Gabber-Kac,Kac} holds, which in this case says that the ideal $D$
is generated by the Serre relations
\begin{align}
({\rm ad}\,e_I)^{1-B_{IJ}}(e_J)=({\rm ad}\,f_I)^{1-B_{IJ}}(f_J)=0\,.
\label{serre0}
\end{align} 
We refer to \cite{Ray} for details about contragredient Borcherds--Kac--Moody superalgebras. 
We also note that different overall normalisations
of the bilinear form $\langle -|-\rangle$ give isomorphic Lie superalgebras $\scr B$ with an isomorphism given by a rescaling 
of $h_0$ and $e_0$. This is however not true for the associated tensor hierarchy algebras below.

Consider the consistent $\mathbb{Z}$-grading of $\scr B$ where $e_0 \in \scr B_1$ and $f_0 \in \scr B_{-1}$, whereas all other generators 
belong to $\scr B_0$. Let $\scr B^\LLoc=\scr B^\LLoc{}_{-1}\oplus \scr B^\LLoc{}_{0} \oplus \scr B^\LLoc{}_{1}$
be the local part of $\scr B$, together with the unique invariant symmetric bilinear form 
such that $\langle x | y \rangle=\kappa(x,y)$ for $x,y\in\fg$.
It then follows from (\ref{matrisdef}) that $\langle e_0 | f_0 \rangle=-\langle f_0 | e_0 \rangle=1$.
Since it is invariant, this form satisfies
\begin{align} \label{nyckelutrakning}
B_{IJ}\langle e_J|f_J\rangle=\langle [h_I,e_J]|f_J\rangle=\langle h_I|[e_J,f_J]\rangle=\langle h_I|h_J\rangle\,,
\end{align}
and since it is gradedly symmetric, this is also equal to $B_{JI}\langle e_I|f_I\rangle $.

Roots are defined for $\scr B$ in the same way as for $\fg$. In particular, the simple roots $\alpha_K$ span the
vector space $\scr H^\ast$ dual to the Cartan subalgebra $\scr H$, which is spanned by the generators $h_K$, where $K=0,1,\ldots,r$.
Let $\varphi:\scr H \to \scr H^\ast$ be the vector space isomorphism given by $\varphi(h_K)=\alpha_K{}^\vee =\langle e_K | f_K\rangle \alpha_K$.
In particular $\alpha_0{}^\vee = \alpha_0$.
It then follows from (\ref{nyckelutrakning}), and the definition $\alpha_J(h_I)=B_{IJ}$ of the simple roots $\alpha_J$,
that we have
\begin{align}
\varphi(h_I)(h_J)= \langle h_I | h_J\rangle\,.
\end{align}
Since $\kappa$ is non-degenerate, $\varphi$ is injective, and $\langle e_K|f_K\rangle \neq0$ for all $K=0,1,\ldots,r$.
We may then introduce an inner product on $\scr H^\ast$ given by
\begin{align}
(\alpha_I,\alpha_J)&=\langle \varphi^{-1}(\alpha_I)|\varphi^{-1}(\alpha_J)\rangle
=\frac1{\langle e_I|f_I\rangle}\frac1{\langle e_J|f_J\rangle}\langle h_I|h_J\rangle=\frac{B_{IJ}}{\langle e_I|f_I\rangle}\,.
\end{align}
In particular,
\begin{align}
(\alpha_0,\alpha_0)&=0\,, & (\alpha_k,\alpha_k)&=\frac2{\langle e_k|f_k\rangle}=\frac2{\kappa(e_k,f_k)}\,
\end{align}
and it follows that
$\alpha_k{}^\vee = 2{\alpha_k}/{(\alpha_k,\alpha_k)}$, as already stated above.

We note that $(\alpha_0{}^\crt,\mu)=-(\lambda,\mu)$ for any $\mu \in \mathfrak{h}^\ast$, since if $\mu=\sum_{i=1}^r m_i\alpha_i$, then
\begin{align}
(\alpha_0{}^\crt,\mu)&=\sum_{i=1}^r m_i(\alpha_0{}^\crt,\alpha_i)=\sum_{i=1}^r m_iB_{0i}=-\sum_{i,j=1}^r \frac{\lambda_i}{\kappa(e_j,f_j)}m_j\delta_{ij}\nn\\
&= -\sum_{i,j=0}^r \frac{\lambda_i}{\kappa(e_j,f_j)}m_j(\Lambda_i,\alpha_j{}^\crt)=
-\sum_{i,j=1}^r {\lambda_i}m_j(\Lambda_i,\alpha_j)=-(\lambda,\mu)\,.
\end{align} 

\begin{prop}
The local part $\scr B^\LLoc$ of $\scr B$ (with respect to the consistent $\mathbb{Z}$-grading) is a contragredient 
local Lie superalgebra where
the element $L$ is given by $L=\sum_{I=0}^r(B^{-1})_{0I}h_I$ and satisfies $\langle L|L\rangle=-1/(\lambda,\lambda)$.
\end{prop}
\Pf \cite{Bossard:2017aae,Cederwall:2017fjm}
We have $[L,e_J]=\alpha_J(L)e_J$, and with $L=\sum_{I=0}^r(B^{-1})_{0I}h_I$ we
get
\begin{align}
\alpha_J(L)=\sum_{I=0}^r(B^{-1})_{0I}\alpha_J(h_I)=\sum_{I=0}^r(B^{-1})_{0I}B_{IJ}=\delta_{0J}\,
\end{align}
as we should. 
Furthermore,
\begin{align}
\langle L | L\rangle &=\sum_{I=0}^r\sum_{J=0}^r(B^{-1})_{0I}(B^{-1})_{0J}\langle h_I | h_J\rangle\nn\\
&=\sum_{I=0}^r\sum_{J=0}^r(B^{-1})_{0I}(B^{-1})_{0J}B_{IJ}\langle e_J | f_J\rangle\nn\\
&=\sum_{J=0}^r\delta_{0J}(B^{-1})_{0J}\langle e_J | f_J\rangle=(B^{-1})_{00}\langle e_0 | f_0\rangle=(B^{-1})_{00}\,.
\end{align}
Since $A$ is invertible,
\begin{align}
(B^{-1})_{00}=\frac{\det{A}}{\det{B}}\neq 0\,,
\end{align}
and in $\mathfrak{h}^\ast$ we can set
\begin{align}
\lambda = \sum_{j=1}^r \frac{(B^{-1})_{j0}}{(B^{-1})_{00}}\alpha_j\,.
\end{align}
We then get
\begin{align}
\lambda_i =
(\alpha_i{}^\vee\!,\lambda)=\sum_{j=1}^r B_{ij}\frac{(B^{-1})_{j0}}{(B^{-1})_{00}}=\sum_{J=0}^r B_{iJ}\frac{(B^{-1})_{J0}}{(B^{-1})_{00}}-B_{i0}
=-B_{i0}\,,
\end{align}
as we should, and
\begin{align}
(\lambda,\lambda)&=\sum_{i=1}^r\sum_{j=1}^r \frac{(B^{-1})_{i0}}{(B^{-1})_{00}}\frac{(B^{-1})_{j0}}{(B^{-1})_{00}}(\alpha_i,\alpha_j)
=\sum_{i=1}^r\sum_{j=1}^r \frac{B_{ij}}{\langle e_i|f_i\rangle}\frac{(B^{-1})_{i0}}{(B^{-1})_{00}}\frac{(B^{-1})_{j0}}{(B^{-1})_{00}}\nn\\
&=-\sum_{i=1}^r \frac{B_{i0}}{\kappa(e_i,f_i)}\frac{(B^{-1})_{i0}}{(B^{-1})_{00}}
=-\sum_{i=1}^r {B_{0i}}\frac{(B^{-1})_{i0}}{(B^{-1})_{00}}\nn\\
&=-\sum_{I=0}^r {B_{0I}}\frac{(B^{-1})_{I0}}{(B^{-1})_{00}}
=-\frac1{(B^{-1})_{00}}=-\frac1{\langle L|L\rangle}\,.
\end{align}
Thus $\langle L|L\rangle=-1/(\lambda,\lambda)$.
\qed

\section{Tensor hierarchy algebras}

In \cite{Carbone:2018xqq}, a Lie superalgebra called {\it tensor hierarchy algebra} and denoted $W$ was associated to any simple and simply laced 
Kac--Moody algebra $\fg$ of rank $r$ and any fundamental weight $\lambda$ of $\fg$. The numbering of fundamental weights of $\fg$
was chosen such that $\lambda=\Lambda_1$. The construction of $W$ started with the 
Cartan matrix $B$ of the Lie superalgebra $\scr B$ associated to the triple $(\fg,\lambda,\kappa)$ as described in the preceding section,
with a normalisation of $\kappa$
such that $\langle e_K|f_K \rangle=1$ for all $K=0,1,\ldots,r$, which means that $B$ is symmetric.
The set of generators $M_{\scr B}=\{e_K,f_K,h_K \,|\, K=0,1,\ldots,r\}$ of $\scr B$ was then modified
to a set $\mangd_{W}$ by replacing the odd generator $f_0$ by $r$ odd generators $f_{0k}$, where $k=0$ or $k=2,3,\ldots,r$.
From this set $\mangd_{W}$ of generators, and the Cartan matrix $B$, an auxiliary Lie superalgebra algebra $\tilde W$
was first constructed as 
the one freely generated by $M_{W}$ modulo the relations
\begin{align}
[h_I,e_J]=B_{IJ}e_J\,,  \qquad   
[h_I,f_J]=-B_{IJ}f_J\,, \qquad   
[e_I,f_J]=\delta_{IJ}h_J\,,  \label{eigen1}
\end{align}\vspace{-18pt}
\begin{align}
({\rm ad}\,e_I)^{1-B_{IJ}}(e_J)=({\rm ad}\,f_I)^{1-B_{IJ}}(f_J)=0\,,
\label{serre1}
\end{align} \vspace{-18pt}
\begin{align}
[e_0,f_{0I}]&=h_I\,, & [h_I,f_{0J}]&=-B_{I0}f_{0J}\,, &
[e_i,[f_j,f_{0K}]]&= \delta_{ij}B_{Kj}f_{0j}\,, \label{eifjf0a}
\end{align} \vspace{-18pt}
\begin{align}
[e_1,f_{0K}]=[f_1,[f_1,f_{0K}]]= 0\,,\label{IdealJ2'}
\end{align}
where $I,J,K=0,1,\ldots,r$ and $i,j,k=0,2,\ldots,r$. (Whenever $f_K$ appears, we assume $K\neq0$, and whenever $f_{0k}$ appears, we assume $k\neq1$.)
Then $W$ was obtained from $\tilde W$ by factoring out the maximal ideal intersecting the local part trivially,
with respect to the consistent $\mathbb{Z}$-grading. By modifying the set of generators further
to $\mangd{}_S=\mangd{}_W\backslash\{h_0,f_{00}\}$ a Lie superalgebra $S$ (called tensor hierarchy algebra as well)
can be defined in the same way (with the relations
involving $h_0$ and $f_{00}$ removed).

It was shown in \cite{Carbone:2018xqq} that $S$
coincides with the original tensor hierarchy algebras introduced in \cite{Palmkvist:2013vya}
in the cases where $\fg$ is finite. It was also shown that $W$ and $S$ coincide with the Lie superalgebras $W(n)$ and $S(n)$ of Cartan type
when $\fg=A_{n-1}$ ($n\geq2$) and $\lambda=\Lambda_1$ with the usual numbering of fundamental weights.
In \cite{Cederwall:2019qnw} the tensor hierarchy algebras were defined in a similar way,
but with the additional relation $[f_{0I},f_{0J}]=0$ in the definition of $\tilde W$, and $W$ obtained from
$\tilde W$ by 
considering a different (non-consistent) $\mathbb{Z}$-grading.

We will now see that the relations (\ref{eigen1})--(\ref{IdealJ2'}) arise naturally in the context of focally associative local algebras.
We consider $\scrB^\LLoc$, the local
part of the contragredient Lie superalgebra
$\scrB$ in the preceding section. We will then investigate the subalgebra of $\scr B^\Loc$ generated by $\scr B_1$ and
$\scr B_{-1}\scr B_0$
modulo the maximal peripheral ideal (where $\scr B_{k}=\scr B^\LLoc{}_{k}$ for $k=\pm1$ and $\scr B_{-1}\scr B_0$ consists
of all products $x_{-1}y_0$ for $x,y\in\scr B^\LLoc$) and 
end the paper with a theorem relating it to
the local part of the tensor hierarchy algebra $W$. 

\begin{prop} \label{apa}
Let $\lambda$ and $\kappa$ be such that {\rm$\lambda^\tri$} is a dominant integral weight. Then
\begin{align}
\big((\alpha_0{}^\crt\!,\alpha)+1\big)\dlb f_0,e_\alpha \drb =0
\end{align}
for all roots $\alpha \neq \alpha_0$ of $\scr B$ with corresponding root vectors $e_\alpha \in \scr B_1$
if and only if $\fg$ is finite and $(\lambda,\theta)=1$ (so that {\rm$\lambda^\tri$} is a pseudo-minuscule weight).
\end{prop}

\Pf Suppose that $\fg$ is infinite-dimensional or that
$(\lambda,\theta)\neq1$.
In either case it is possible to find a root $\zeta$ 
of $\fg$ 
such that $\alpha_0+\zeta$ is a root of $\scr B$ and $(\alpha_0{}^\crt,\zeta)\neq-1$,
for example $\zeta=\theta$ if $\fg$ is finite, since then
\begin{align}
(\alpha_0{}^\crt,\zeta)=(\alpha_0{}^\crt,\theta)
=-(\lambda,\theta)\neq -1\,.
\end{align}
We can then set $\alpha=\alpha_0+\zeta$ and $e_\alpha=\dlb e_0,e_\zeta \drb$, where $e_\zeta$ is a root vector
corresponding to $\zeta$, and it follows that $\big((\alpha_0{}^\crt,\alpha)+1\big)\dlb f_0,e_\alpha \drb \neq0$.

On the other hand, suppose that $\fg$ is finite and that
$\lambda^\tri$ is a pseudo-minuscule weight, say $\lambda^\tri=\Lambda_j$ for some $j$ such that $c_j=1$
and let $\alpha$ be a root of $\scr B$ such that $e_\alpha \in \scr B_1$ and $\dlb f_0,e_\alpha \drb \neq0$.
Then $\alpha-\alpha_0=\sum_{k=1}^r b_k \alpha_k$ is a root of $\fg$ and
\begin{align}
(\alpha_0{}^\crt,\alpha)=(\alpha_0{}^\crt,\alpha-\alpha_0)=\sum_{i=1}^rb_i(\alpha_0{}^\crt,\alpha_i)=-\sum_{i=1}^rb_i\lambda^\tri{}_i
=-\sum_{i=1}^rb_i\delta_{ij}=-b_j\,.
\end{align}
Since $\theta$ is the highest root, $b_j \leq c_j=1$. But $b_j$ cannot be zero since then $(\alpha_0{}^\crt,\alpha-\alpha_0)$ would be zero
as well, and $\alpha=\alpha_0+(\alpha-\alpha_0)$ would not be a root. Thus $b_j=1$ and $(\alpha_0{}^\crt,\alpha)+1=0$.
\qed
\\
\noindent
In what remains of this paper, we assume that $\fg$ is a finite Kac--Moody algebra.

\begin{lemma} \label{sistalemmat}
Let $\lambda$ and $\kappa$ be such that $\lambda^\tri$ is a pseudo-minuscule weight.
Then $f_0(h_0+L)$
generates a peripheral ideal of the subalgebra generated by $\scr B_1$ and $\scr B_{-1}\scr B_0$.
\end{lemma}
\Pf
Let $e_\alpha \in \scr B_1$ be a root vector of a root $\alpha \neq \alpha_0$. We then have
\begin{align}
[e_0,f_0 h_0]+[e_0,f_0 L]&= [e_0,f_0]h_0-f_0[e_0,h_0]+[e_0,f_0]L-f_0[e_0,L]\nn\\
&=h_0+L+f_0e_0=h_0+L-h_0-L=0\,,\nn\\[8pt]
[e_\alpha,f_0h_0]+[e_\alpha,f_0L]&=[e_\alpha,f_0]h_0-f_0[e_\alpha,h_0]+[e_\alpha,f_0]L-f_0[e_\alpha,L]\nn\\&=f_0[h_0,e_\alpha]+f_0[L,e_\alpha]\nn\\
&=\big((\alpha_0{}^\vee\!,\alpha)+1\big)f_0e_\alpha=-\big((\alpha_0{}^\vee\!,\alpha)+1\big)\dlb f_0,e_\alpha\drb=0\,,
\end{align}
where the last equation follows from Proposition \ref{apa}.
\qed

\begin{theorem} \label{satsen}
Let $\lambda$ and $\kappa$ be such that $\lambda^\tri$ is a pseudo-minuscule weight, and choose a numbering of the fundamental weights
of $\fg$ such that $\lambda^\tri=\Lambda_1$.
Then there is a local Lie superalgebra homomorphism
from the local part of $W$ to the subalgebra of 
 $\scr B^\Loc$ generated by $\scr B_1$ and $\scr B_{-1}\scr B_0$
modulo the maximal peripheral ideal,
given by $f_{0K} \mapsto f_0 h_K$ (and leaving the other generators unchanged).
\end{theorem}

\Pf Let $V$ be the subalgebra of $\scr B^\Loc$ generated by $\scr B_1$ and $\scr B_{-1}\scr B_0$,
and let 
$D$ be the maximal peripheral ideal of $V$.
We will show that the relations (\ref{eigen1})--(\ref{IdealJ2'}) are satisfied in $V/D$ with $f_{0K}=f_0 h_K$.
We will first
show that $[f_0h,e_1]=[f_1,[f_1,f_0h]]=0$ for any $h \in \mathfrak{h}$.
From Lemma \ref{sistalemmat} we know that $f_0(h_0+L)=0$ in $V/D$, and we also have $B_{01}=-1$ since $\lambda^\tri=\Lambda_1$.
We get
\begin{align}
[f_0h,e_1]&=f_0[h,e_1]=\alpha_1(h)f_0e_1\nn\\
&=-{\alpha_1(h)}f_0[h_0+L,e_1]=-{\alpha_1(h)}[f_0(h_0+L),e_1]=0\,.
\end{align}
Similarly,
\begin{align} \label{410}
[f_1,[f_1,f_0h]]&=[f_1,f_0[f_1,h]]+[f_1,[f_1,f_0]h]\nn\\
&=f_0[f_1,[f_1,h]]+2[f_1,f_0][f_1,h]+[f_1,[f_1,f_0]]h\nn\\
&=2[f_1,f_0][f_1,h]=2\alpha_1(h)[f_1,f_0]f_1\nn\\
&=-2{\alpha_1(h)}[f_1,f_0][f_1,h_0+L]\,.
\end{align}
We can then perform the first three steps in (\ref{410}) backwards, but with $h$ replaced by $h_0+L$,
and find that
\begin{align}
[f_1,f_0][f_1,h_0+L]=
[f_1,[f_1,f_0(h_0+L)]]=0\,. 
\end{align}
We have thus shown the relations (\ref{IdealJ2'}).
The other relations involving $f_{0K}$ are straightforward to show,
\begin{align}
[e_0,f_{0K}]&=[e_0,f_{0}h_{K}]=[e_0,f_0]h_K-f_0[e_0,h_K]=h_K\,,\nn\\
[h_I,f_{0J}]&=[h_I,f_{0}h_{J}]=[h_I,f_{0}]h_J+f_0[h_I,h_J]=-B_{I0}f_0h_J=-B_{I0}f_{0J}\,,\nn\\
[e_i,[f_j,f_{0K}]]&=[e_i,[f_j,f_{0}h_{K}]]
=f_0[e_i,[f_j,h_K]]=\de_{ij}B_{Kj}f_0h_j=\de_{ij}B_{Kj}f_{0j}\,,
\end{align}
and those not involving $f_{0K}$ automatically satisfied.
\qed

\noindent
We conjecture that this homomorphism is in fact an isomorphism in the simply laced case,
but leave the proof for future work.
\\\\
\noindent
\underline{\it Acknowledgments:} JP would like to thank Per Bäck, Jens Fjelstad and Victor Kac for discussions and 
helpful answers to questions.

\bibliographystyle{utphysmod2}

%\bibliography{biblio.bib}

%\bibliography{biblio}

\providecommand{\href}[2]{#2}\begingroup  %\raggedright

%\endgroup

\end{document}